\documentclass{article}
\pdfoutput=1	
\usepackage{fullpage}	
\usepackage[breaklinks,hidelinks]{hyperref}
\hypersetup{
  pdftitle={Idealizing Useful Fictions in Omega Grounded Arithmetic},
  pdfauthor={Bryan Ford}}

\newcommand{\oga}{OGA\xspace}
\newcommand{\ati}{\textsc{ati}\xspace}
\newcommand{\PT}[1]{\mathrm{PT}(#1)}
\providecommand{\godelnum}[1]{\ulcorner #1 \urcorner}

\usepackage{amsmath}	
\usepackage{amssymb}	
\usepackage{amsthm}

\newtheorem{thm}{Theorem}[section]
\newtheorem{cor}{Corollary}[section]
\newtheorem{lem}{Lemma}[section]
\newtheorem{rem}{Remark}[section]

\usepackage{times}
\usepackage{xspace}
\usepackage{cite}		
\usepackage{bussproofs}
\usepackage{stmaryrd}
\usepackage[noabbrev,capitalize]{cleveref}
\usepackage{thmtools}
\usepackage{xstring}
\usepackage{graphicx}
\usepackage{subcaption}
\usepackage{nicematrix}		
\usepackage{longtable}
\usepackage{pifont}

\usepackage{isabelle}
\usepackage{isabellesym}
\usepackage{pdfpages}
\usepackage{ragged2e}
\newcommand{\SNIP}[2]{\expandafter\newcommand\csname snippet--#1\endcsname{#2}}
\IfFileExists{snips.tex}{\input{snips}}{}

\newcommand{\GetSnip}[1]{%
    \ifcsname snippet--#1\endcsname%
        \csname snippet--#1\endcsname%
    \else%
        \PackageWarning{snips}{Snippet ``#1'' is undefined.}%
        \emph{Warning: Snippet ``#1'' is undefined.}%
    \fi%
}

\newcommand{\Snippet}[1]{{%
  \newcount\i
  \i=0
  \loop
    \GetSnip{#1-\the\i}%
    \advance \i 1
  \ifcsname snippet--#1-\the\i\endcsname
  \repeat
}}

\newcommand{\SnippetPart}[3]{{%
  \newcount\i
  \i=#1
  \loop
    \ifnum \i=#2
      \renewcommand{\isanewline}{}%
    \fi
    \GetSnip{#3-\the\i}%
    \advance \i 1
    \ifnum \i>#2 {}
    \else \repeat
}}

\hypersetup{breaklinks=true}

\newcommand{\com}[1]{}

\long\def\note#1{}	

\newcommand{\pga}{PGA\xspace}		
\newcommand{\rga}{RGA\xspace}		
\newcommand{\pa}{PA\xspace}		

\newcommand{\limp}{\rightarrow}		
\newcommand{\liff}{\leftrightarrow}	
\newcommand{\ldef}{\equiv}		

\newcommand{\typestyle}[1]{\textsf{#1}\xspace}
\com{	
\newcommand{\tbool}{\typestyle{bool}}
\newcommand{\ttrue}{\typestyle{true}}
\newcommand{\tfalse}{\typestyle{false}}

\newcommand{\tnat}{\typestyle{nat}}

}
\newcommand{\tbool}{\typestyle{B}}	
\newcommand{\ttrue}{\typestyle{T}}
\newcommand{\tfalse}{\typestyle{F}}

\newcommand{\tnat}{\typestyle{N}}

\newcommand{\ctrue}{\ttrue}
\newcommand{\cfalse}{\tfalse}

\newcommand{\judgment}[1]{\textsf{ #1}}

\newcommand{\jbool}{\judgment{\tbool}}

\newcommand{\jnat}{\judgment{\tnat}}

\makeatletter
\newlength{\doublefracgap}
\DeclareRobustCommand{\doublefrac}[2]{%
  \mathinner{\mathpalette\doublefrac@{{#1}{#2}}}%
}
\newcommand{\doublefrac@}[2]{\doublefrac@@#1#2}
\newcommand{\doublefrac@@}[3]{%
  \ooalign{%
    \raisebox{\doublefracgap}{$\m@th#1\frac{#2}{\phantom{#3}}$}\cr
    \raisebox{-\doublefracgap}{$\m@th#1\frac{\phantom{#2}}{#3}$}\cr
  }%
}
\newcommand{\ddoublefrac}[2]{{\displaystyle\doublefrac{#1}{#2}}}

\makeatother

\newcommand{\irl}[1]{\ensuremath{\mathit{#1}}}		

\newcommand{\infrule}[3][]{\cfrac{#2}{#3}\IfStrEq{#1}{}{}{\ \irl{#1}}}
\newcommand{\infeqv}[3][]{\ddoublefrac{#2}{#3}\IfStrEq{#1}{}{}{\ \irl{#1}}}
\newcommand{\infceqv}[4][]{\cfrac{#2\qquad}{}\ddoublefrac{#3}{#4}
				\IfStrEq{#1}{}{}{\ \irl{#1}}}

\newcommand{\z}{\mathbf 0}

\newcommand{\suc}{\mathbf S}

\newcommand{\cR}{\mathcal R}

\newcommand{\tofn}[1]{\underline{#1}}

\crefname{thm}{Theorem}{Theorems}
\crefname{cor}{Corollary}{Corollaries}
\crefname{lem}{Lemma}{Lemmas}
\crefname{prop}{Proposition}{Propositions}
\crefname{defi}{Definition}{Definitions}
\crefname{rem}{Remark}{Remarks}
\crefname{exa}{Example}{Examples}

\newcommand{\cA}{\mathcal A}		
\newcommand{\hs}{\mathrel{\triangleright}}	
\newcommand{\lft}[1]{{\uparrow}#1}	
\newcommand{\vz}{v_0}			
\newcommand{\wf}{\mathrm{wf}}		
\newcommand{\fv}{\mathrm{fv}}		

\begin{document}

\title{Idealizing Useful Fictions in Omega Grounded Arithmetic}
\author{Bryan Ford \\ EPFL}
\date{}					
\maketitle
Grounded arithmetic is a family of formal systems for reasoning
about computation in which a statement may be asserted only when a
terminating computation backs it; the logics are \emph{paracomplete}
--- for a sentence whose backing computation never settles, neither
the sentence nor its negation is derivable, so paradoxes like the
Liar are harmless rather than explosive.  The reflective member of
the family, RGA, can quantify over its own computations, but leaves
one thing conspicuously open: it cannot certify that its own
unbounded searches have definite yes-or-no answers.  This paper
studies what happens when that openness is closed by exactly one
rule --- ATI, the $\omega$-grounded universal: if every numeric
instance of a universal sentence is certified decided, the universal
is certified decided.  The resulting system, OGA, shares RGA's
syntax and rules symbol-for-symbol otherwise, and every consequence
is developed as a machine-checked theorem.  Decidedness certificates
become abundant --- every totality question about a computable
function is certified to have an answer, whether or not anyone can
produce it --- and this is exactly the provable separation between
the two systems.  The rule costs no completeness: OGA is complete
for its own semantics, in which certified-but-unresolved sentences
receive values built from the system's own open questions.
Provability remains recursively enumerable, with a
primitive-recursive certificate checker, while $\omega$-truth
deliberately is not --- and within that asymmetry, incompleteness
takes a new form.  The G\"odel sentence is classified,
unconditionally, as a \emph{genuine fiction}: neither provable nor
refutable, yet valued, and carrying a computable \emph{pedigree}
recording exactly what adopting it as an axiom commits one to.  The
adoption is then itself a theorem suite: extending OGA by any
finite stock of true fictions is consistent, and independently
certified adoptions can never collide.  All results are formalized
in Isabelle/HOL with no unproven assumptions.

%
%

\section{Introduction}
\label{sec:intro}

Formal reasoning about computation has an old irritant: classical
logic demands that every sentence be true or false, but a
computation is entitled to a third option --- running forever.
Grounded arithmetic~\cite{ford24reasoning,ford25have} takes the
third option seriously.  It replaces the classical picture of truth
--- every sentence true or false in advance, provability chasing
that pre-existing totality of facts --- with an earned one: a
sentence holds when a terminating semantic process \emph{grounds}
it, and a sentence whose backing computation never settles simply
holds no value, harmlessly.  Reflective Grounded Arithmetic
(RGA)~\cite{rga-paper} showed that this discipline can support
genuine quantification over infinite domains, by making the
quantifiers computational: a universal sentence is represented by a
\emph{search sentence} --- a formula that runs the corresponding
unbounded search, stage by stage, over a compiled step function ---
so that the system's own searches become ordinary objects of its
arithmetic.  (\emph{Reflective} refers to exactly this: the system
reasons about its own computations.)

RGA leaves one thing conspicuously open.  Its search sentences are
\emph{undecided}: the system cannot prove, of its own searches,
that they have a definite yes-or-no answer (in the family's
vocabulary, that they are \emph{decided}, carrying a
\emph{decidedness certificate}).  This is not an accident of
presentation but a theorem --- a search that never succeeds has no
terminating computation backing any verdict about it, and
certifying such a search as decided would collapse RGA's
discipline that provability means groundedness~\cite{rga-paper}.  The
question this paper answers is what arithmetic looks like when that
openness is closed deliberately, by the weakest natural means: a
single rule.

\paragraph{One rule.}  The rule --- $\omega$-grounded universal
decidedness introduction, ATI --- says: if every numeric instance
of a universal sentence's body is decided, then the universal
sentence itself is decided.  One schematic derivation covering the
instances suffices; and the same premise serves existential
sentences too, since deciding an existential, no less than a
universal, requires surveying every instance.  The resulting system, OGA,
shares RGA's syntax symbol-for-symbol and its proof system rule-for-rule,
except for ATI.  We present the two systems together as a
delta (\cref{sec:system}), so the reader sees the entire difference
at a glance --- and then watches everything else in the paper fall
out of it.

\paragraph{Not a trade.}  It would be natural to expect the new rule
to be bought with something --- for a system that decides more to be
complete for less.  It is not.  RGA is complete for its own grounded
semantics; OGA is complete for \emph{its} own semantics
(\cref{sec:semantics}), in which a certified sentence whose answer
is out of reach receives a \emph{value} built from the system's
own open questions rather than a truth-value gap --- and
provability again coincides exactly with the semantics' top
verdict.  This is the family
pattern rather than a coincidence: each grounded system is complete
for the semantics obtained by reflecting its own characteristic
infinitary move, and adding a rule moves the system and its semantics
together.  What changes between RGA and OGA is not how much of its
semantics the system can reach, but which semantics it has ---
the same discipline over a richer reflected ground.  The one thing
that does come apart is $\omega$-truth, and that separation is
deliberate and confined (\cref{sec:metatheory}).

\paragraph{What falls out.}  First, the searches become decided
(\cref{sec:decidedness}): OGA proves of its own search sentences
that they have answers --- instance by instance, and then for every
totality question (``does this computable function halt on every
input?'')\ outright, with no hypothesis on the function.  This is
precisely the separation between the systems: measured against a
generic ladder of quantifier strengths shared by the whole family,
OGA occupies the rung that RGA provably cannot reach.  Second, the metatheory splits along a designed
asymmetry (\cref{sec:metatheory}): OGA-provability is recursively
enumerable, with a primitive-recursive judgment-level certificate
checker, while $\omega$-truth is deliberately \emph{not} r.e.
Third, within that asymmetry incompleteness itself is reclassified.
The G\"odel sentence of OGA is not merely ``true but unprovable'':
it is, unconditionally and by machine-checked proof, a
\emph{genuine fiction} --- neither provable nor refutable, yet
carrying a definite value in the fact/fiction semantics of
\cref{sec:semantics}, together with a computable \emph{pedigree}
recording exactly which open questions adopting it as an axiom
would rest on.  The
classical trichotomy of provable/refutable/undecidable becomes a
semantic classification the system itself can reason about.

\paragraph{The price, if any.}  Decidedness may yet not be free, but
the bill --- if it comes --- is not completeness.  We present
evidence (\cref{sec:concl}) that the classical L\"ob obstruction
--- the theorem that blocks a strong-enough classical system from
proving its own consistency --- is structurally inert in RGA
precisely because provability there is not a decided matter, yet
reactivates in OGA: the stronger system plausibly loses
self-verification exactly by gaining decidedness.  Resolving that
inversion, and the graded set theory that OGA's certified searches
support, are the subjects of companion work in progress.

\paragraph{Contributions.}
\begin{itemize}
\item The system OGA: RGA plus the single ATI rule, presented as a
  delta, with design rationale (\cref{sec:system}).
\item A supervaluation $\omega$-semantics with soundness, a
  certificate-gated fafi semantics for which OGA is \emph{complete},
  and a ladder of derived decidedness principles topped by the
  unconditional decidedness of every totality question
  (\cref{sec:semantics}).
\item The separation theorem: decidedness of the reflective ground
  characterizes the RGA/OGA gap, at generic-tier level
  (\cref{sec:decidedness}).
\item The metatheoretic asymmetry: r.e.\ provability with a
  primitive-recursive checker versus non-r.e.\ $\omega$-truth, and
  $\omega$-incompleteness (\cref{sec:metatheory}).
\item The unconditional classification of the G\"odel sentence as a
  genuine fiction with computable pedigree (\cref{sec:metatheory}).
\item The adoption calculus for fictions: extending \oga{} by any
  finite stock of $\omega$-true sentences is consistent, certified
  adoptions merge without conflict (confluence), and the G\"odel
  sentence assembles into a single adoptable-fiction bundle ---
  decided, irrefutable, unprovable, fictional, safe to adopt
  (\cref{sec:metatheory:adopt}).
\item A complete Isabelle/HOL formalization with no unproven
  assumptions (\cref{sec:formal}).
\end{itemize}

\textbf{AI disclosure:} both the mechanically-verified proofs and
the writing of this paper were significantly assisted by artificial
intelligence (Claude Fable from Anthropic), as detailed in
\cref{sec:formal:ai}.

\section{Background: Grounded Arithmetic and \rga}
\label{sec:bg}

\emph{Grounded arithmetic} replaces the classical stipulation that
every arithmetic sentence is true or false with a constructive
discipline: a sentence holds when a terminating semantic process
grounds it, fails when such a process refutes it, and otherwise has
no truth value at all.  Provability is designed to track groundedness
rather than to approximate a presumed-total truth predicate.  The
program develops this discipline as a ladder of systems of
increasing quantificational strength; \rga~\cite{rga-paper} is the
tier at which genuine universal quantification first appears, and
the present paper's system \oga{} is its one-rule extension.

\paragraph{\rga{} in brief.}
\rga{} grounds universal quantification \emph{reflectively}: the
system reasons about its own computations.  A computation is given
as a \emph{step function}, a total function reporting the state of
a staged process at each step.  Every primitive-recursive step
function is compiled to a term of the arithmetic itself (the
compiler \texttt{cmp}).  A universal sentence is then represented
by a \emph{search sentence} over the compiled function: a formula
that runs the corresponding unbounded search, stage by stage.  The
universal counts true when the system's own proof search certifies
its schematic instance.  It counts false when the system refutes a
particular numeral instance, in which case an explicit
counterexample witness is extractable.  We call the stock of
representable search sentences the system's \emph{reflective
ground}.  It makes \rga's own searches first-class objects of
arithmetic, and it is the raw material for everything in this
paper.

\paragraph{Metatheorems of \rga{} (statements only).}
The following are proved, machine-checked, in~\cite{rga-paper};
we use them here as black boxes.

\begin{itemize}
\item \emph{Soundness} (\texttt{rga\_sound}): derivable judgments
  are grounded true under every assignment satisfying their
  hypotheses.  \emph{Consistency} (\texttt{rga\_consistent})
  follows: no sentence is both provable and refutable.
\item \emph{Open completeness / the truth--provability coincidence}
  (\texttt{rga\_complete\_open}, \texttt{rga\_open\_iff}): for
  well-formed sentences, open or closed, provability coincides with
  grounded truth.  Grounded truth is not an ideal that provability
  approximates.  The two are the same relation.
\item \emph{N-soundness} (\texttt{rga\_N\_sound}): every provable
  totality claim is backed by an actual value at an actual step
  index.
\item \emph{Expressive power} (\texttt{rga\_truth\_re},
  representation): grounded truth is \emph{recursively enumerable}
  --- some program lists exactly the true sentences --- and the
  sets of numbers \rga{} can represent are exactly the recursively
  enumerable ones.  This is the Church--Turing landmark: the
  system's reach coincides with computability itself.
\item \emph{$\omega$-incompleteness}: some family of statements
  has every numeric instance provable while its universal closure
  is not merely unprovable but semantically ungrounded.
\item \emph{Substructural character}.  Double-negation elimination
  holds.  Excluded middle for quantified sentences fails.  A
  refuted universal always yields a counterexample witness (a
  Markov-flavored corner of the system).  And the deduction theorem
  --- the classical principle that from ``$\Gamma$ plus $a$ proves
  $c$'' one may conclude ``$\Gamma$ proves $a$ implies $c$'' ---
  fails in exactly one place: hypotheses $a$ that are not decided.
  This last point matters repeatedly below.  In grounded systems,
  implication introduction demands a \emph{decided} antecedent.
\end{itemize}

\paragraph{The undecided reflective ground.}
One further fact about \rga{} frames this paper.  Call a sentence
\emph{decided} when the system proves it has a definite yes-or-no
answer (the judgment $\jbool$, the \emph{decidedness certificate}
of \cref{sec:intro}).  \rga's search sentences are, in general,
\emph{not} provably decided.  This is a theorem, not an open
question: a search that never settles grounds no verdict, so
certifying such a search as decided would collapse \rga's
discipline that provability means groundedness
(\texttt{rga\_rrep\_not\_decided}%
\footnote{Proved in the accompanying formalization, postdating the
conference version of~\cite{rga-paper}.  
keep here with this footnote, or restate in \cref{sec:decidedness}
as part of the separation theorem?
}).
\oga{} is what results from closing exactly this openness,
deliberately, by the weakest natural means.

\section{The Systems \rga{} and \oga{}, Together}
\label{sec:system}

\oga{} shares its syntax with \rga{} symbol for symbol, and its
proof system rule for rule --- except for exactly one additional
rule.\footnote{One bookkeeping caveat, for readers comparing
against the formalizations.  \rga's mechanized rule set carries
three additional legacy primitives (equation forms of
$\beta$-reduction and the recursor's zero case) that are derivable
from the conversion rules.  \oga's rule set omits them; they
remain admissible via the embedding of \rga{} derivations.  The
forty-two rules shown in the figures are the shared set both
systems interderive.}  We therefore present both systems at once.
\Cref{fig:syntax} gives the common syntax;
\cref{fig:rules,fig:rules-reflect} give the forty-two shared rules;
and \cref{fig:ati} gives the one rule that separates the systems.
A reader familiar with~\cite{rga-paper} may skip directly
to \cref{sec:system:ati}.

\subsection{Shared Syntax and Judgments}
\label{sec:system:shared}

\begin{figure}
\centering
$
t \ldef v
  \mid \bot
  \mid \z
  \mid \suc t
  \mid t = t
  \mid \neg t
  \mid t \lor t
  \mid \lambda t
  \mid t \cdot t
  \mid \cR
  \mid \cA
$
\medskip

\begin{tabular}{r@{~$\ldef$~}l@{\qquad}r@{~$\ldef$~}l}
$t \ne u$ & $\neg(t = u)$ &
$t \jnat$ & $t = t$ \\
$t \jbool$ & $t \lor \neg t$ &
$t \land u$ & $\neg(\neg t \lor \neg u)$ \\
$t \limp u$ & $\neg t \lor u$ &
$t \liff u$ & $(t \limp u) \land (u \limp t)$ \\
$\forall.\,b$ & $\cA \cdot (\lambda b)$ &
$\exists.\,b$ & $\neg(\cA \cdot (\lambda \neg b))$ \\
\end{tabular}
\caption{Term syntax of \rga{} and \oga{} --- identical (top) and derived forms (bottom).
  Variables $v$ are de Bruijn indices and $\lambda$ is the only
  binder; the quantifiers bind nothing themselves, so the derived
  $\forall.\,b$ quantifies the de Bruijn variable $0$ of the
  body $b$.}
\label{fig:syntax}
\end{figure}

Terms comprise variables (de Bruijn indices --- variables are
numbered rather than named), the divergent constant $\bot$, zero
and successor, the grounded connectives (equality, negation,
disjunction), $\lambda$-abstraction with application, and two
combinator constants: the primitive recursor $\cR$ and the
universal quantifier $\cA$.  One caution for readers arriving from
other traditions: $\bot$ is not a falsity proposition.  It is a
\emph{term} whose evaluation never terminates, the canonical
divergent computation.  The derived forms --- the type judgments
$t \jnat$ and $t \jbool$, the remaining connectives, and the
quantifiers --- are abbreviations, exactly as in \rga.

Judgments take the form $\Gamma \vdash c$ with $\Gamma$ a finite
\emph{set} of hypothesis terms.  A \emph{proof} is a finite list of
judgments, each following from earlier entries by a single rule; the
judgment $\Gamma \vdash c$ is \emph{derivable} when some valid proof
contains it.  Proofs are finite syntactic objects and hence
encodable as numbers (G\"odel-codable).  Both the reflective
semantics and the certificate checkers of \cref{sec:metatheory}
quantify over exactly these codes.  The rule names and groupings
below match the formalization, where the two systems' rule sets are
kept textually aligned case by case so that metatheoretic inductions
transfer between them.

\begin{figure}
\centering
\textit{structural rules}\\[0.5ex]
$\infrule[hyp]{\wf(c)}{\Gamma, c \vdash c}$
\qquad
$\infrule[wk]{\Gamma \vdash c}{\Gamma, a \vdash c}$
\qquad
$\infrule[cut]{\Gamma \vdash a \quad \Gamma, a \vdash c}
  {\Gamma \vdash c}$
\qquad
$\infrule[sub]{\Gamma \vdash c \quad \wf(\Gamma,c,s)}
  {[v{\mapsto}s]\Gamma \vdash [v{\mapsto}s]c}$
\medskip

\textit{truth values are numerals}\\[0.5ex]
$\infeqv[T{=}]{\Gamma \vdash a}{\Gamma \vdash a = \tofn{1}}$
\qquad
$\infeqv[F{=}]{\Gamma \vdash \neg a}{\Gamma \vdash a = \tofn{0}}$
\medskip

\textit{negation}\\[0.5ex]
$\infrule[{\neg}E]{\Gamma \vdash \neg a \quad \Gamma \vdash a
    \quad \wf(c)}
  {\Gamma \vdash c}$
\qquad
$\infeqv[{\neg}{\neg}]{\Gamma \vdash a}{\Gamma \vdash \neg\neg a}$
\medskip

\textit{disjunction}\\[0.5ex]
$\infrule[{\lor}I_1]{\Gamma \vdash a \quad \wf(b)}
  {\Gamma \vdash a \lor b}$
\qquad
$\infrule[{\lor}I_2]{\Gamma \vdash b \quad \wf(a)}
  {\Gamma \vdash a \lor b}$
\qquad
$\infrule[{\lor}I_3]{\Gamma \vdash \neg a \quad \Gamma \vdash \neg b}
  {\Gamma \vdash \neg(a \lor b)}$
\\[1ex]
$\infrule[{\lor}E_1]{\Gamma \vdash a \lor b \quad
    \Gamma, a \vdash c \quad \Gamma, b \vdash c}
  {\Gamma \vdash c}$
\qquad
$\infrule[{\lor}E_2]{\Gamma \vdash \neg(a \lor b)}
  {\Gamma \vdash \neg a}$
\qquad
$\infrule[{\lor}E_3]{\Gamma \vdash \neg(a \lor b)}
  {\Gamma \vdash \neg b}$
\medskip

\textit{grounded equality}\\[0.5ex]
$\infeqv[{=}refl]{\Gamma \vdash a \jnat}{\Gamma \vdash a = a}$
\qquad
$\infrule[{=}sym]{\Gamma \vdash a = b}{\Gamma \vdash b = a}$
\qquad
$\infrule[{=}subs]{\Gamma \vdash a = b \quad
    \Gamma \vdash [v{\mapsto}a]p \quad \wf(p)}
  {\Gamma \vdash [v{\mapsto}b]p}$
\\[1ex]
$\infrule[{\ne}sym]{\Gamma \vdash a \ne b}{\Gamma \vdash b \ne a}$
\qquad
$\infrule[{=}TI]{\Gamma \vdash a \jnat \quad \Gamma \vdash b \jnat}
  {\Gamma \vdash (a = b) \jbool}$
\qquad
$\infrule[{=}TE]{\Gamma \vdash (a = b) \jbool}
  {\Gamma \vdash (a \jnat) \land (b \jnat)}$
\medskip

\textit{natural numbers and induction}\\[0.5ex]
$\infrule[\z{}I]{}{\Gamma \vdash \z \jnat}$
\qquad
$\infrule[\suc{}TI]{\Gamma \vdash a \jnat}
  {\Gamma \vdash \suc a \jnat}$
\qquad
$\infrule[\suc{}TE]{\Gamma \vdash \suc a \jnat}
  {\Gamma \vdash a \jnat}$
\qquad
$\infrule[\suc{}{\ne}\z]{\Gamma \vdash a \jnat}
  {\Gamma \vdash \suc a \ne \z}$
\\[1ex]
$\infrule[\suc{}{\ne}I]{\Gamma \vdash a \ne b}
  {\Gamma \vdash \suc a \ne \suc b}$
\qquad
$\infrule[\suc{}{=}E]{\Gamma \vdash \suc a = \suc b}
  {\Gamma \vdash a = b}$
\qquad
$\infrule[\suc{}{\ne}E]{\Gamma \vdash \suc a \ne \suc b}
  {\Gamma \vdash a \ne b}$
\qquad
$\infrule[\bot{}TE]{\Gamma \vdash \bot \jnat \quad \wf(p)}
  {\Gamma \vdash p}$
\\[1ex]
$\infrule[ind]{\Gamma \vdash [i{\mapsto}\z]c \quad
    v_i \jnat,\, c,\, \Gamma \vdash [i{\mapsto}\suc v_i]c \quad
    \Gamma \vdash a \jnat \quad
    \wf(c) \quad i \notin \fv(\Gamma)}
  {\Gamma \vdash [i{\mapsto}a]c}$
\caption{The \rga proof rules, part one:
  the structural and logical trunk shared with \pga.
  A double line means the rule holds in both directions
  (each direction is one primitive rule);
  $\wf(\cdot)$ marks well-formedness side conditions,
  which attach only to terms not already
  occurring in a premise judgment.}
\label{fig:rules}
\end{figure}

\begin{figure}
\centering
\textit{syntactic head steps}\\[0.5ex]
$\infrule[\beta]{\wf(t)}{(\lambda b) \cdot t \hs b[t]}$
\qquad
$\infrule[\cR\z]{\wf(h)}{\cR \cdot g \cdot h \cdot \z \hs g}$
\qquad
$\infrule[\cR\suc]{\wf(g)}
  {\cR \cdot g \cdot h \cdot (\suc \tofn{n}) \hs
   h \cdot \tofn{n} \cdot (\cR \cdot g \cdot h \cdot \tofn{n})}$
\qquad
$\infrule[spine]{t \hs t'}{t \cdot u \hs t' \cdot u}$
\medskip

\textit{conversion along head steps}\\[0.5ex]
$\infrule[conv]{\Gamma \vdash c' \quad c \hs c'}
  {\Gamma \vdash c}$
\qquad
$\infrule[convE]{\Gamma \vdash c \quad c \hs c'}
  {\Gamma \vdash c'}$
\\[1ex]
$\infrule[{\neg}convE]{\Gamma \vdash \neg c \quad c \hs c'}
  {\Gamma \vdash \neg c'}$
\qquad
$\infrule[conv{=}]{\Gamma \vdash c' \jnat \quad c \hs c'}
  {\Gamma \vdash c = c'}$
\medskip

\textit{schematic recursor unfolding}\\[0.5ex]
$\infrule[\cR\suc{}{=}]{\Gamma \vdash m \jnat \quad
    \Gamma \vdash (h \cdot m \cdot (\cR \cdot g \cdot h \cdot m)) \jnat
    \quad \wf(g)}
  {\Gamma \vdash \cR \cdot g \cdot h \cdot (\suc m) =
    h \cdot m \cdot (\cR \cdot g \cdot h \cdot m)}$
\medskip

\textit{the reflective quantifier}\\[0.5ex]
$\infrule[{\cA}I]{\vz \jnat,\, \lft\Gamma \vdash \lft f \cdot \vz}
  {\Gamma \vdash \cA \cdot f}$
\qquad
$\infrule[{\cA}E]{\Gamma \vdash \cA \cdot f \quad \Gamma \vdash a \jnat}
  {\Gamma \vdash f \cdot a}$
\\[1ex]
$\infrule[{\neg\cA}I]{\Gamma \vdash \neg(f \cdot a) \quad
    \Gamma \vdash a \jnat}
  {\Gamma \vdash \neg(\cA \cdot f)}$
\qquad
$\infrule[{\neg\cA}E]{\Gamma \vdash \neg(\cA \cdot f) \quad
    \vz \jnat,\, \neg(\lft f \cdot \vz),\, \lft\Gamma \vdash \lft c}
  {\Gamma \vdash c}$
\caption{The \rga proof rules, part two:
  computation and the universal quantifier.
  The head-step relation $t \hs t'$ is a
  primitive-recursively decidable syntactic relation,
  not a judgment;
  $b[t]$ is the $\beta$-instance substituting $t$ for
  the outermost bound variable of $b$,
  $\tofn{n}$ is a syntactic numeral,
  and $\lft{(\cdot)}$ shifts all free de Bruijn indices up by one,
  so that $\vz$ is fresh for $\lft\Gamma$, $\lft f$, and $\lft c$.}
\label{fig:rules-reflect}
\end{figure}

\subsection{The One New Rule: \ati}
\label{sec:system:ati}

\begin{figure}
\centering
$\infrule[ATI]{\Gamma \vdash \cA \cdot
    (\lambda\, ((\lft u \cdot \vz) \jbool))}
  {\Gamma \vdash (\cA \cdot u) \jbool}$
\caption{The single rule separating \oga{} from \rga:
  $\omega$-grounded universal decidedness introduction (\ati).
  The premise is a universally quantified decidedness claim about
  the body $u$; the conclusion decides the quantified sentence
  itself.}
\label{fig:ati}
\end{figure}

The rule \ati{} (\cref{fig:ati}) says: if every numeric instance of
the body is decided --- stated as an object-logic universal over the
$\jbool$-judgment of the body --- then the universally quantified
sentence is itself decided.  Three design points deserve emphasis.

\paragraph{The premise is internal.}  \ati{} is a finitary rule.
Its premise is a single \oga{}-sentence --- an internal universal
over decidedness claims --- not an infinitary scheme of
meta-premises.  \oga{} does not adopt the $\omega$-rule; it adopts
the reflection of one specific $\omega$-fact, decidedness, into the
object logic.

\paragraph{One premise serves both quantifiers.}  The
existential case --- deciding $\exists.\,b$ --- is derivable from
\ati{} with the \emph{same} universal premise.  The derived
existential is a negated universal, and decidedness is closed
under negation.  Deciding an existential, no less than a
universal, requires surveying every instance.  There is no cheaper
warrant: in $\omega$-grounded quantifier semantics, a quantified
sentence's status is determined by its instances alone, with no
additional oracle to consult.

\paragraph{Decidedness, not truth.}  \ati{} concludes
$(\cA \cdot u) \jbool$, never $\cA \cdot u$.  The rule adds no new
grounded truths at all on its own; what it adds is Boolean
\emph{status}.  Every consequence developed in this paper ---
the separation from \rga, the metatheoretic asymmetry, the
reclassification of incompleteness --- flows from status, not from
new truths.  This is also why \oga's consistency and
$\omega$-soundness replay \rga's case for case
(\cref{sec:semantics}): the one new case preserves exactly the
invariants the other forty-two already maintained.

\paragraph{Excluded middle, without a witness.}  Unfolding the
abbreviation of \cref{fig:syntax} makes the rule's character plain.
$(\cA \cdot u) \jbool$ \emph{is}
$(\cA \cdot u) \lor \neg(\cA \cdot u)$, so what \ati{} concludes is
literally an instance of excluded middle for the quantified
sentence.  \oga{} is thus \rga{} plus classical excluded middle at
the quantifier --- admitted not unrestrictedly, but under a
certificate: the premise demands that every instance already be
decided.

The reading is worth making explicit because it identifies what
kind of principle \ati{} is.  A constructivist objection to
excluded middle is that asserting $\varphi \lor \neg \varphi$
without possessing either disjunct claims a completed survey one
has not performed.  \ati{} answers the objection in the only way a
grounded system can: it demands the survey, in the form of the
pointwise premise, and grants the disjunction on the strength of
it.  What it still does not supply is a witness.  The semantics of
\cref{sec:semantics} makes this exact rather than figurative ---
the disjunction evaluates to $\top$ while neither disjunct does
(\cref{lem:totality-lem}) --- so \oga{} models
reasoning-without-witnesses as an object-level phenomenon, with
the price of each such assertion recorded in the pedigree of the
sentence it licenses.

\section{Semantics: Facts, Fictions, and the $\omega$-Model}
\label{sec:semantics}

\oga{} has two semantics, playing different roles.  The
\emph{$\omega$-completion} extends \rga's grounded semantics with
the one clause \ati{} reflects.  Its job is safety: it delivers
soundness and consistency.  The \emph{fafi semantics}
(``fact/fiction'') is \oga's own native semantics, and it delivers
the paper's conceptual vocabulary.  It is a
\emph{supervaluation}: a scheme in which a sentence containing
unresolved unknowns counts as true when \emph{every} admissible
way of resolving those unknowns would make it true.  Ours is
moreover \emph{certificate-gated}: an unknown may enter the
picture only when the system has derived a certificate for it.
Under this semantics every valued sentence is either a \emph{Fact}
or a \emph{Fict}, and every Fict carries a computable
\emph{pedigree} recording exactly what adopting it rests on.

\subsection{The $\omega$-completion and soundness}
\label{sec:semantics:omega}

\begin{thm}[$\omega$-soundness; \texttt{omega\_sound}]
\label{thm:omega-sound}
If $\Gamma \vdash c$ in \oga{} and every hypothesis in $\Gamma$ is
$\omega$-true, then $c$ is $\omega$-true.
\end{thm}
Consistency follows as in \rga.  The one new proof case --- \ati{}
--- is exactly the statement that supervaluation respects
$\omega$-grounded decidedness, previewing the machinery below.

\subsection{The fafi value algebra}
\label{sec:semantics:fafi}

Values live in a Boolean algebra generated by \emph{presumption
atoms}: named unknowns, one per certified-but-unresolved question.
A value is a Boolean combination of the constants \ctrue, \cfalse
and these atoms.  Grounded sentences evaluate to the constants
(\emph{Facts}).  A universally quantified sentence whose body is
certified decided --- and only such a sentence --- may introduce
the named unknown $x_u$ as its value (a \emph{Fict}).  The gate is
a proof: the atom for $\cA \cdot u$ enters only with an
\oga-derivation of the pointwise decidedness premise, mirroring
\ati{} exactly.  This gating is the reflective discipline applied
a second time.  \rga{} applied it to the $\omega$-clause for
quantifiers; \oga{} applies it to supervaluation itself.  The
gating also keeps the semantics computable in reach: classical
supervaluation quantifies over arbitrary completions of the
language, far beyond anything a program could survey, while the
gated form stays recursively enumerable.  One more property
matters throughout, and readers who know three-valued logic will
recognize what is at stake.  In a truth-functional scheme such as
Kleene's, \emph{unknown or unknown} is unknown --- even when the
two unknowns are the same question, asked positively and
negatively.  Here the two occurrences stay correlated:
$(\cA \cdot u) \lor \neg(\cA \cdot u)$ computes to
$x_u \lor \neg x_u = \ctrue$, a Fact.  That computation is exactly
\ati's soundness obligation, and \cref{lem:totality-lem} records
the case of it that the later sections use.

\subsection{The coincidence: provable = super-true}
\label{sec:semantics:c1}

\begin{thm}[fafi characterisation;
  \texttt{c1\_characterisation}]
\label{thm:c1}
For well-formed $\varphi$:
$\emptyset \vdash \varphi$ in \oga{} iff $\varphi$ is
\emph{super-true} --- true under every admissible valuation of the
atoms, where a valuation is admissible when it never contradicts a
verdict the system has actually derived.
\end{thm}

\oga{} is thus sound and \emph{complete} for its own
reflective-supervaluational truth.  This extends the family
pattern: each grounded system is complete for the semantics
obtained by reflecting its own characteristic infinitary move.
The proof has two load-bearing ideas.  One is the exact definition
of admissibility just glossed.  The other is a maximality argument
(via Zorn's lemma) over consistent sets, engineered to sidestep
the deduction theorem's decided-antecedent restriction.  Both are
sketched in \cref{sec:semantics:appx}.
The constants fragment coincides with \rga's reflective semantics:
the forty-two shared rules compute Facts.

\subsection{The trichotomy and the pedigree}
\label{sec:semantics:fict}

Every well-formed sentence of \oga{} falls into exactly one of
three classes:
\begin{itemize}
\item \emph{Fact}: a constant value --- grounded, and provable or
  refutable accordingly;
\item \emph{Fict}: a committed indeterminate --- valued, neither
  provable nor refutable, introduced through the decidedness gate;
\item \emph{no value}: the gate never fires.  The Liar is the
  paradigm (\texttt{feval\_liar\_gap}).  \oga{} remains a gap
  logic: supervaluation does not paper over ungroundedness.
  Self-reference is not the only source, though.  The gate fires
  on universals only, so an ungrounded sentence of any other shape
  --- a redex, for instance, a term with a pending computation
  step --- gaps as well (\cref{sec:decidedness:fict}).
\end{itemize}

The third class is not a residue, and the line it draws is sharper
than ``ungrounded''.  Two ungrounded sentences can land in
different classes, and which one they land in tracks something
real.  Take a diagonal gap point: the universal
$\cA \cdot \mathit{ofam}\,x$ is ungrounded, but its pointwise
decidedness is provable, so the gate fires unconditionally
(\texttt{ofam\_decided}) and the sentence receives an atom.  It is
ungrounded and yet fictionally decided.  The Liar is ungrounded
too, and receives nothing at all: no valuation whatever assigns it
a value (\texttt{feval\_liar\_gap}).

The difference is not an artifact of where the gate happens to
fire.  At the diagonal point there is an answer and only a
certificate is missing --- which is exactly why
\cref{sec:metatheory:fict} can say that we know that answer from
outside while \oga{} does not.  At the Liar there is no answer to
be missing, and nothing attaches because there is nothing for it
to attach to.  So \ati{} admits a one-line characterization: it
converts absence of a \emph{certificate} into a fiction, and
leaves absence of a \emph{fact} exactly where it found it.

One qualification, since the conversion is sensitive to how a
sentence is presented as well as to what it says.  A search-sentence
instance at a diagonal point is a missing-certificate case too, but
it is a redex rather than a universal, so the gate does not reach
it and it gaps --- one $\beta$-step away from a value
(\cref{sec:decidedness:fict}).  Gaps of presentation and gaps of
content both sit in the third class; only the second kind is beyond
\ati's reach in principle.

A value's \emph{support} is finite and consists of decided
universal atoms (\texttt{feval\_support}); its \emph{pedigree} is
the Fict part of that support --- the sentence's ideal commitments.
A Fict atom commits exactly to itself
(\texttt{pedigree\_fatom\_fict}); a Fact commits to nothing
(\texttt{pedigree\_fatom\_prov}).  Pedigree makes proof-theoretic
provenance \emph{semantic}: to adopt a fiction is to accept a
finite, computable set of ideal commitments, and the algebra
tracks them compositionally.  This is the formal content of the
paper's title: \oga's fictions are \emph{useful} because they are
valued, consistent to adopt, and priced.

\subsection{The decidedness ladder}
\label{sec:semantics:ladder}

\ati{} is stated for a single quantified sentence, but the work of
the following sections is done by a short derived ladder, each rung
obtained from the one below by a fixed manoeuvre.  We record it
here; \cref{sec:decidedness,sec:metatheory} consume it throughout.

The bottom rung is pointwise.  A representing term applied to a
numeral head-steps to the negation of a universal whose body is a
quantifier-free equation on a compiled term.  That body is decided
at every instance by the schematic premise machinery, so \ati{}
applies, and the decidedness transports back across the
conversion.

\begin{lem}[$\Sigma_1$ decidedness; \texttt{rrep\_decided}]
\label{lem:rrep-decided}
For every primitive-recursive step function $f$ and every $n$:
$\emptyset \vdash (\mathit{rrep}\,f \cdot \tofn{n}) \jbool$.
\end{lem}

The same argument never inspects the numeral.  Running it at an
arbitrary term known to be a number, in an arbitrary context, gives
the schematic rung --- the form that survives passage under a
binder:

\begin{lem}[schematic $\Sigma_1$ decidedness;
  \texttt{rrep\_decided\_gen}]
\label{lem:rrep-decided-gen}
If $\Gamma \vdash Y \jnat$ then
$\Gamma \vdash (\mathit{rrep}\,f \cdot Y) \jbool$.
\end{lem}

Instantiating $Y$ at a fresh variable of type $\tnat$, discharging
it by universal introduction, and applying \ati{} once more lifts
the ladder over the quantifier itself:

\begin{thm}[totality is decided; \texttt{totality\_decided}]
\label{thm:totality-decided}
For every primitive-recursive step function $f$ ---
\emph{whether or not its search is total} ---
$\emptyset \vdash (\cA \cdot \mathit{rrep}\,f) \jbool$.
\end{thm}

\Cref{thm:totality-decided} is the ladder's working top, and the
form most of this paper actually uses: \oga{} settles the
\emph{status} of every totality question, unconditionally in $f$.
It claims no truth --- \ati's conclusion is always a decidedness
judgment, never an assertion of the universal --- and that gap
between deciding and affirming is exactly what
\cref{sec:metatheory} is about.  Two riders come free by the same
route: the normalization guard of the compiled fragment is decided
(\texttt{normT\_decided}), as is the diagonal family
(\texttt{cmp\_diag\_decided}), which is what makes the diagonal
arguments of \cref{sec:metatheory} unconditional.

Decidedness is not merely a status label: it is the key that
unlocks implication introduction.  Grounded systems grant the
deduction theorem only for \emph{decided} antecedents
(\cref{sec:bg}), so every rung of this ladder is simultaneously a
licence to discharge a hypothesis.  One mechanism, several
consequences: it is what fires the fafi certificate gate above, and
what makes the tier-2 construction of \cref{sec:metatheory:tier2}
go through where the corresponding \rga{} construction stalls.  We
flag it at each appearance.


\section{Decidedness of the Reflective Ground}
\label{sec:decidedness}

This section proves the paper's headline structural claim: the one
rule \ati{} decides \emph{every} search sentence --- and this is
exactly the boundary between \oga{} and \rga, in the strong sense
that \rga{} provably cannot cross it.

\subsection{The Separation Theorem}
\label{sec:decidedness:sep}

The decidedness ladder of \cref{sec:semantics:ladder} already
carries the positive half.  Its bottom rung
(\cref{lem:rrep-decided}) decides every search sentence at every
point, and its top (\cref{thm:totality-decided}) decides every
totality question outright, unconditionally in the step function.
What this section adds is the other half: that \rga{} can have
neither.

\begin{thm}[undecidedness at \rga;
  \texttt{rga\_rrep\_not\_decided}, from~\cite{rga-paper}'s
  formalization]
\label{thm:rrep-not-decided}
There is an $n$ with
$\emptyset \nvdash (\mathit{rrep}\,\mathit{sfApp} \cdot \tofn{n})
\jbool$ in \rga.
\end{thm}

The witness is the halting diagonal: for $x$ outside the halting
set, the search sentence at $\langle x,x\rangle$ is semantically
ungrounded, and \rga's provability never outruns groundedness.

\begin{cor}[separation]
\label{cor:separation}
\oga's theorems strictly extend \rga's: the sentence of
\cref{thm:rrep-not-decided} is \oga-provable and not
\rga-provable.  The extension consists of \emph{status}
judgments: \ati{} adds no grounded truths
(\cref{sec:system:ati}), yet its status judgments are new theorems.
\end{cor}

\subsection{The Tier Reading}
\label{sec:decidedness:tier}

The formalization maintains a generic ladder of quantifier tiers
--- computation, propositional, predicate, and
$\omega$-decidedness --- as \emph{locales}: named bundles of
assumed inference rules, abstract interfaces independent of any
particular system, which a concrete system \emph{realizes} by
proving the bundled rules.  The $\omega$-decidedness tier is characterized
by the signature pair: universal and existential decidedness
introduction from the \emph{same} pointwise premise.  \oga{}
registers at this tier (the formalization's \texttt{OGA\_Tier}):
both interface obligations are discharged from \ati{} through a
conversion bridge.  \Cref{thm:rrep-not-decided} says \rga{}
provably cannot register: any system with \rga's
truth--provability coincidence that decided its own searches would
certify an ungrounded sentence.  Decidedness of the reflective
ground is thus not an incremental strengthening but a tier
boundary, and \ati{} is exactly the rule that crosses it.

\subsection{Status Without Truth: Fiction at the Universal}
\label{sec:decidedness:fict}

The two readings of the separation meet in the fafi semantics.
They meet one level up from the separation witness itself,
however, and it is worth being exact about the distance between
the two levels.

\paragraph{The universal is a fiction.}
Consider a totality universal $\cA \cdot \mathit{rrep}\,f$.  It is
decided, for every $f$ without exception, by
\cref{thm:totality-decided}.  Decidedness is precisely what fires
the certificate gate of \cref{sec:semantics:fafi}.  So the gate
fires here, and the sentence receives a value:

\begin{lem}[totality universals are valued;
  \texttt{feval\_totality\_fict}]
\label{lem:totality-fict}
For every $f$, the sentence $\cA \cdot \mathit{rrep}\,f$ takes as
a value the indeterminate it names: writing $u = \mathit{rrep}\,f$,
the atom $x_{\cA \cdot u}$ is among its values.
\end{lem}

The atom fires for every $f$ without exception --- but that by
itself does not make the sentence a fiction, and the distinction
is worth drawing carefully.  Evaluation is a \emph{relation}, not
a function.  Where the sentence is also grounded --- $f$ a step
function whose totality \oga{} can actually prove --- it carries a
constant value as well, and is a Fact; the atom's availability
costs nothing and settles nothing.  A \emph{genuine} Fict is the
stricter case in which the atom fires and \emph{nothing grounds
the sentence at all}.

That stricter case is what the diagonal supplies, and it is worth
being clear about which ingredient does the work.  Decidedness is
free, from \ati, for every $f$ alike.  Unprovability is not: it is
the Gödel content, imported at the diagonal
(\cref{sec:metatheory:fict}), and it is what turns an available
indeterminate into a committed one.  In that case, and only there,
the pedigree is the sentence's own singleton --- adopting it
commits one to exactly itself and to nothing further.

For search universals the stricter case admits an exact
recursion-theoretic description, which we record because it says
precisely where \oga's fictions live:

\begin{thm}[fictions sit in the gap; \texttt{fict\_iff\_gap}]
\label{thm:fict-iff-gap}
The sentence $\cA \cdot \mathit{rmat}\,f\,n$ --- that the search
of $f$ at $n$ never hits --- is a genuine fiction exactly when the
search really never hits, and \rga{} does not prove that it never
hits.
\end{thm}

The first condition is a $\Pi_1$ truth --- a statement of the
form ``for every $n$, a mechanically checkable fact holds.''  The
second is the failure of a recursively enumerable approximation to
reach it.  So the fictions of the reflective ground are precisely
the points of the $\Pi_1$-minus-r.e.\ difference: universal
truths that no listing of theorems attains.  That difference is a semantic object, and it is worth
saying so, since our witnesses for it are all diagonal: being a
diagonal point is how one \emph{exhibits} a member of the
difference, not what \emph{makes} one.  A gap point carries no
record of how it was found.  This is also worth stating carefully
for what the theorem does
\emph{not} say.  It relates fictionality to a gap between
\emph{objects} --- a set of numbers and a proof system's reach ---
and says nothing whatever about how any theorem about those
objects was proved.  Fact/Fict status is a property of the formula
and the calculus, invariant under the metatheorist's choice of
argument; a sentence does not become fictional because we reasoned
classically about it.

The excluded middle over it, meanwhile, is a Fact outright:

\begin{lem}[excluded middle is recovered;
  \texttt{feval\_totality\_lem}]
\label{lem:totality-lem}
For every $f$, the sentence
$(\cA \cdot \mathit{rrep}\,f) \lor \neg(\cA \cdot \mathit{rrep}\,f)$
evaluates to \ctrue.
\end{lem}

This is the disjunction property failing, made semantic: the
disjunction is \ctrue{} while neither disjunct is.  Correlation
between the two occurrences of the atom is what buys it, and it is
exactly what truth-functional three-valued schemes cannot deliver.
It is also the precise form of the observation in
\cref{sec:system:ati}: $\jbool$ abbreviates a disjunction with a
negation, so \cref{lem:totality-lem} is excluded middle for the
sentence holding outright while the sentence itself remains
unsettled --- the assertion of a classical disjunction with no
witness for either side, valued rather than merely stipulated.

\paragraph{The instance is not.}
The pointwise separation witness of \cref{thm:rrep-not-decided}
behaves quite differently, and the difference is instructive.
\oga{} decides that witness too, by \cref{lem:rrep-decided}.  Yet
it may have no value at all.

The reason is that decidedness and valuedness are different
achievements.  Decidedness is a judgment of the proof system.  A
fafi value, by contrast, is computed by induction on syntactic
shape --- and the certificate gate is one clause of that
induction, the clause for universals.  A representing term applied
to a numeral is not a universal.  It is a redex.  The gate
therefore has nothing to offer it, and such an instance is valued
exactly when it is \emph{grounded}.  At a diagonal point it is not
grounded, and so it gaps.

One consequence is worth recording, since it is easy to assume
otherwise: \oga's gaps are not confined to self-reference.  The
Liar (\cref{sec:semantics:fict}) is the paradigm case of a gap,
but it is not the only one.

\paragraph{The gap is one step deep.}
The instance gaps, but only just.  Contracting the redex yields
$\neg(\cA \cdot \mathit{rmat}\,f\,n)$, and that sentence is valued
unconditionally:

\begin{lem}[the contractum is always valued;
  \texttt{rmat\_univ\_feval}]
\label{lem:reduct-valued}
For every $f$ and every $n$, the sentence
$\neg(\cA \cdot \mathit{rmat}\,f\,n)$ evaluates to
$\neg x_{\cA \cdot \mathit{rmat}\,f\,n}$ --- whether or not the
instance it came from has a value.
\end{lem}

What a gapping instance lacks, then, is not content that the
semantics is unable to supply.  It is a shape that the gate is
able to recognize.  The value sits one $\beta$-step away, and the
evaluation relation does not take that step: it is defined by
induction on five syntactic forms, with no conversion clause.  The
asymmetry is between presentation and content.

It is worth naming because the proof system does not share it.
There, conversion transports decidedness freely --- indeed
\cref{lem:rrep-decided} is itself proved through this very head
step, reaching the instance from the contractum that
\cref{lem:reduct-valued} evaluates.  The same step that the
calculus takes for granted is the one the semantics declines.

\medskip\noindent
What \oga{} gains over \rga{} at the tier boundary is therefore
sharper than ``its searches become ideal elements.''  Its
unresolved \emph{totality questions} become ideal elements:
valued, adoptable, priced.  \Cref{sec:metatheory} shows that
G\"odel's incompleteness, in this setting, is an instance of
exactly that phenomenon.

\section{Metatheory: Enumerability, Incompleteness, and Fictional
  Status}
\label{sec:metatheory}

\rga's metatheory closed with a coincidence: grounded truth and
provability are one r.e.\ relation.  \oga{} deliberately breaks the
coincidence in one direction only.  Provability remains r.e., with
an explicitly primitive-recursive certificate checker; semantic
truth strictly outruns it.  Incompleteness then returns --- but
with a different face: the classical ``true but unprovable''
sentence reappears as a sentence whose \emph{correct}
classification, provable inside the semantics, is \emph{adoptable
fiction}.

\subsection{Provability is Recursively Enumerable}
\label{sec:metatheory:re}

\begin{thm}[\texttt{oga\_entails\_re}, \texttt{ogaPfCheck}]
\label{thm:re}
\oga-derivability of coded judgments is recursively enumerable:
there is a primitive-recursive certificate test accepting exactly
the codes of derivable judgments, with proof codes as
certificates.
\end{thm}

The checker follows the one-step architecture used throughout the
formalization: a proof is a list of judgments, each entry justified
from its suffix by one rule, and the per-entry test is primitive
recursive.  \ati{} adds one case to the checker exactly as it adds
one rule to the calculus.  The checker is exposed at two
granularities: a per-judgment test (\texttt{sfPfJ}) and a search
wrapper over whole proofs (\texttt{sfPfO}).  Both are consumed
downstream --- by the premise machinery of \cref{sec:decidedness}
and, twice more, by the companion set-theory development.

Enumerability has an internal face.  Let $\mathit{Prov}$ be the
search sentence over the proof-checking step function --- \oga's own
provability predicate, an ordinary term of the arithmetic.

\begin{thm}[the provability predicate is exactly right;
  \texttt{oga\_prov\_characterises}]
\label{thm:prov-char}
$\emptyset \vdash \mathit{Prov}(\godelnum{\varphi})$ \;iff\;
$\emptyset \vdash \varphi$.
\end{thm}

The forward direction is the \emph{admissibility of the provability
rule}: whenever \oga{} proves that $\varphi$ is provable, $\varphi$
is provable --- a reflection-flavored principle obtained without
internal reflection, since the argument runs through
$\omega$-soundness at the meta level: a provable search sentence
really hits, and the hit \emph{is} a proof.  The backward direction
is the first Hilbert--Bernays derivability condition (D1, in the
standard numbering: a proof of $\varphi$ can always be converted
into a proof that $\varphi$ is provable).  The forward direction
is D1's converse, which \pa{} enjoys only as a consequence of
$\Sigma_1$-soundness and cannot itself prove.  No L\"ob-style
obstruction applies to admissibility, which is a statement about
the meta-level rule and not about an internal
implication.\footnote{We claim nothing here about the remaining
derivability conditions: neither the internal distribution
$\mathit{Prov}(\godelnum{\varphi \limp \psi}) \limp
(\mathit{Prov}(\godelnum{\varphi}) \limp
\mathit{Prov}(\godelnum{\psi}))$ nor internal
$\Sigma_1$-completeness, in the form
$\mathit{Prov}(\godelnum{\varphi}) \limp
\mathit{Prov}(\godelnum{\mathit{Prov}(\godelnum{\varphi})})$,
is established in the formalization, and nothing below depends on
them.  Whether they hold --- and hence whether L\"ob's theorem
itself is available in \oga{} --- is the open question of
\cref{sec:concl}.  One constraint is already visible: were both
available, L\"ob's theorem would apply, and an internal soundness
schema $\mathit{Prov}(\godelnum{\varphi}) \limp \varphi$ would
then make \oga{} inconsistent.  Internal soundness cannot be had
together with the other two conditions.}

\subsection{Decidedness of the Provability Predicate}
\label{sec:metatheory:decd}

$\mathit{Prov}$ is a search sentence like any other, so the ladder
of \cref{sec:semantics:ladder} applies to it.  This yields a
property no classical provability predicate has in the same way:

\begin{cor}[provability is decided; \texttt{prov\_decided},
  \texttt{frset\_mem\_decided}]
\label{cor:prov-decided}
For every $\varphi$: $\emptyset \vdash
\mathit{Prov}(\godelnum{\varphi}) \jbool$.  More generally the
membership formula of any r.e.\ set is decided, uniformly in both
the set and the element.
\end{cor}

In a classical theory this is trivial, and worthless: excluded
middle gives it for every formula whatever.  Here it is neither.
\oga{} is a gap logic in which decidedness must be earned, the
Liar does not have it, and the reflective ground does not have it
pointwise (\cref{sec:decidedness:fict}).  That $\mathit{Prov}$ has
it, unconditionally and at every instance, is a substantive fact
about the shape of the provability predicate: it is a search
sentence, and \ati{} decides those.  Note also what it is
\emph{not}: a decided sentence need not be an affirmed one, so
\cref{cor:prov-decided} settles the status of
$\mathit{Prov}(\godelnum{\varphi})$ without settling its truth ---
it is a Fict, in the vocabulary of \cref{sec:semantics:fict}, not
a Fact.

The property is inherited by the compounds one wants to reason
with.  Write
\[
  \mathit{Decd}(\godelnum{\varphi}) \;\ldef\;
  \mathit{Prov}(\godelnum{\varphi}) \lor
  \mathit{Prov}(\godelnum{\neg\varphi})
\]
for the claim that $\varphi$ is settled one way or the other.
Since $\jbool$ is closed under disjunction, two copies of
\cref{cor:prov-decided} give:

\begin{cor}[decidedness is decided; \texttt{fdecd\_decided}]
\label{cor:decd-decided}
For every $\varphi$: $\emptyset \vdash
\mathit{Decd}(\godelnum{\varphi}) \jbool$.
\end{cor}

So \oga{} settles, for every sentence without exception, the
question \emph{``is this sentence settled?''}  We take this up in
\cref{sec:metatheory:revenge}, where it turns out to bear on a
familiar difficulty for gappy theories of truth.

One direction of internalization is nonetheless blocked, and it is
worth recording alongside the positive results, since it marks the
edge of what \cref{thm:prov-char} delivers:

\begin{thm}[completeness does not internalize;
  \texttt{prov\_complete\_not\_internal}]
\label{thm:prov-not-internal}
There is a $\varphi$ with $\emptyset \vdash \varphi \jbool$ and
$\emptyset \nvdash \varphi \limp
\mathit{Prov}(\godelnum{\varphi})$.
\end{thm}

The witness is the halting diagonal once more: there $\varphi$ is
$\omega$-true while $\mathit{Prov}(\godelnum{\varphi})$ is
$\omega$-false, so the implication is $\omega$-false and
$\omega$-soundness forbids proving it.  Restricting to decided
$\varphi$ does not rescue it --- the witness is decided, by
\cref{thm:totality-decided}.  The two halves of
\cref{thm:prov-char} are therefore not symmetric under
internalization: the admissibility half has a free positive
branch, while the completeness half has an outright
counterexample.  D1 holds of \oga, but D1 does not internalize
\emph{in} \oga.

\subsection{Revenge, and Where \oga{} Does and Does Not Escape It}
\label{sec:metatheory:revenge}

Gappy theories of truth are haunted by \emph{revenge}.  Kripke's
least fixed point~\cite{kripke75outline} leaves the Liar without a
truth value, which is the intended outcome; but the theory cannot
say so in its own language, and a language that could express
``ungrounded'' would generate a Liar the fixed point does not
classify.  Field's response~\cite{field08saving} adds a
determinacy operator and pays for it with a transfinite hierarchy
of them, each level classifying the sentences the level below left
open.  The general shape is familiar: the vocabulary used to
\emph{describe} the gap re-enacts the paradox, and the repair is
an ascent~\cite{beall08revenge}.

\Cref{cor:decd-decided} looks like it refuses the ascent.  \oga{}
settles ``is $\varphi$ settled?'' for every $\varphi$, and it does
so unconditionally --- the corollary is a schema in $\varphi$ with
a uniform one-line proof.  In particular it applies when $\varphi$
itself mentions $\mathit{Decd}$.  So
$\mathit{Decd}(\godelnum{\mathit{Decd}(\godelnum{\psi})})$ is
decided, by the same argument that decides
$\mathit{Decd}(\godelnum{\psi})$, and so on upward.  Where Field
has a hierarchy that does not collapse, \oga{} has a family of
statements that are all instances of one lemma.  The hierarchy is
flat.

We think the flatness is real but that it must be attributed
correctly, because it is easy to claim too much for it.  \oga{}
has not defused the Liar: the Liar still gaps
(\cref{sec:semantics:fict}), and \oga{} remains a paracomplete
theory in the same business as its predecessors.  What is going on
is rather that $\mathit{Decd}$ \emph{is not the gap predicate}.
$\mathit{Decd}$ reports on the proof system.  It is a disjunction
of search sentences --- a $\Sigma_1$ object, one asserting that
some checkable witness exists --- and \ati{} decides exactly such
objects.  The predicate
that would play Field's role is a predicate of \emph{semantic
value}: one saying that a sentence has a fafi value at all rather
than gapping.  \oga{} does not internalize that predicate.  Its
evaluation relation is a meta-level construction, defined by
induction outside the object language, and nothing in the calculus
expresses ``this sentence has no value.''

So \oga's position is intermediate rather than victorious, and the
intermediate position is the interesting one.  Kripke's theory
cannot internally express groundedness at all.  Field's can
express determinacy and pays in transfinite ascent.  \oga{}
internally expresses a decidedness predicate which is
\emph{total on the reflective ground} and flat under iteration,
while declining to internalize semantic valuedness --- the one
predicate whose internalization is what makes revenge bite.  The
Tarskian retreat is still there.  It has been moved, and the
region it has been moved out of is exactly the region where the
classical theories needed their hierarchies.

Two qualifications keep this honest.  First, decided is not
affirmed.  \Cref{cor:decd-decided} licenses reasoning by cases on
whether $\varphi$ is settled --- this is the deduction-theorem
unlock of \cref{sec:semantics:ladder}, and it is a genuine
inferential gain --- but \oga{} does not in general prove either
branch.  The trichotomy is available as a reasoning device, not as
a decision procedure, and \cref{thm:prov-not-internal} shows the
gap between the two is not an artifact of presentation.  Second,
the revenge test proper has not been run.  Revenge needs a
determinacy predicate that is both internal and subject to
diagonalization, and $\mathit{Decd}$ as defined here is a
meta-level operation on $\varphi$ rather than an internally
computed function of $\godelnum{\varphi}$.  Making it the latter
--- which the diagonal lemma would require --- is not
difficult in principle, but it has not been done, and until it is,
\cref{cor:decd-decided} should be read as evidence that the ascent
is unnecessary rather than as proof that it is impossible.

\subsection{$\omega$-Truth is Not --- by Design}
\label{sec:metatheory:asym}

In \rga, grounded truth is r.e.\ (\texttt{rga\_truth\_re}) and
coincides with provability.  \oga's $\omega$-truth cannot be r.e.:
it decides every search sentence (\cref{lem:rrep-decided} is sound
for it), so an enumeration of $\omega$-truth would enumerate the
complement of the halting diagonal --- a set no program can list.  The asymmetry is the design:
the \emph{calculus} stays inside the enumerable world --- every
theorem carries a checkable certificate --- while the
\emph{semantics} certifies status one level beyond it.  The gap
between the two is not a defect to be closed but the exact home of
the G\"odel phenomena, to which we now turn.

\subsection{$\omega$-Incompleteness, Sharpened}
\label{sec:metatheory:incomplete}

\begin{thm}[\texttt{oga\_omega\_incomplete}]
\label{thm:incomplete}
There is a closed, well-formed $u$ with
$\emptyset \vdash u \cdot \tofn{n}$ for every $n$, while
$\emptyset \nvdash \cA \cdot u$.
\end{thm}

The witness family is again the halting diagonal
(\texttt{oga\_ouniv\_gap}, \texttt{rrep\_univ\_gap}).  The
sharpening over \rga's $\omega$-incompleteness is the status of
the unprovable closure: in \rga{} it was semantically
\emph{ungrounded}; in \oga{} it is \emph{decided}
(\cref{thm:totality-decided}) --- the system proves it has a
Boolean answer --- and still neither provable nor refutable.  Instancewise
provability, certified Boolean status, no verdict: this is
precisely the profile the fafi semantics names.

\subsection{Termination Provability}
\label{sec:metatheory:pt}

The incompleteness family has a quantitative face.  Define the
\emph{termination-provability classes}: $\PT{S}$ is the set of
primitive-recursive step functions whose search \emph{everywhere
hits} provably in $S$ --- $f \in \PT{S}$ iff
$\emptyset \vdash_S \cA \cdot \mathit{rrep}\,f$.

\begin{thm}[\texttt{PT\_incl}, \texttt{PT\_sound},
  \texttt{PT\_ceiling}, \texttt{PT\_ceiling\_rga}]
\label{thm:pt}
$\PT{\rga} \subseteq \PT{\oga}$; every member of $\PT{\oga}$
really hits everywhere; and there is an everywhere-hitting step
function outside $\PT{\oga}$ --- hence also outside $\PT{\rga}$.
\end{thm}

The ceiling witness is the halting diagonal again: at a nonhalting
point the search hits at every argument --- at step zero, even ---
yet no \oga-proof affirms the universal.  Termination provability
thus stratifies strictly below actual termination, and it does so
on both sides of the inclusion: the ceiling is a property of
r.e.\ affirmation as such, not an artifact of \ati.

That symmetry must not be misread as progress toward separation.
The \rga-side ceiling is obtained by \emph{contraposing} the
embedding --- every \rga{} derivation is an \oga{} derivation, so
an \oga{} non-theorem is a fortiori an \rga{} non-theorem --- and
the witness is therefore the same function in both statements.  It
exhibits a point where \rga{} fails \emph{because \oga{} already
fails there}; what strictness would need is a point where \rga{}
fails and \oga{} succeeds, and we have none.  The inclusion
$\PT{\rga} \subseteq \PT{\oga}$ is rule containment and we claim no
more: whether it is \emph{strict} is open.  \ati's conclusion is
always a decidedness judgment and never an affirmation, so a
separating function would have to exploit decidedness indirectly
--- through the unlocked deduction theorem, as in
\cref{sec:metatheory:tier2} --- and we have no candidate.  We
deliberately frame this analysis without ordinals: in classical
systems, termination provability and ordinal analysis go hand in
hand, but whether that correspondence survives in grounded systems
is open, and nothing here presupposes it.

\subsection{The G\"odel Sentence is a Genuine Fiction}
\label{sec:metatheory:fict}

\begin{thm}[\texttt{oga\_godel\_fict},
  \texttt{oga\_totality\_fict\_uncond}]
\label{thm:godel-fict}
There is a diagonal instance $x$ whose universal sentence
satisfies $\mathit{fict}$: it is valued by the fafi semantics,
neither provable nor refutable, and its pedigree is its own
singleton.  In particular some totality claim
$\cA \cdot \mathit{rrep}\,f$ is a genuine fiction.
\end{thm}

The valuedness half is \cref{lem:totality-fict}, and holds of every
totality universal whatever; what the diagonal supplies is the
other half, that this particular one is settled in neither
direction.

Three features distinguish this from the classical statement of
incompleteness.  \emph{Unconditionality}: no consistency or
$\omega$-consistency hypothesis appears --- the diagonal argument
runs against the halting set and the landed soundness theorems,
and the classification is absolute.  \emph{Status, not truth}:
the theorem does not say the G\"odel sentence is ``true but
unprovable'' --- $\omega$-truth of the instance family is a
metatheoretic fact, while the object-level classification is
\emph{valued, undecided, adoptable}.  \emph{The price is
computable}: the pedigree (\cref{sec:semantics:fict}) of the
G\"odel sentence is its own singleton --- adopting it as an axiom
commits one to exactly that sentence and nothing else.
Incompleteness, in the grounded-supervaluational setting, is the
statement that the system manufactures its own ideal elements:
sentences whose correct epistemic treatment --- certified by the
semantics --- is adoption as priced fictions rather than assertion
as truths.

\paragraph{What the fiction is about.}
The second of those features repays being made concrete, because
it is easy to state backwards.  Two vantage points are in play at
once, and $\emptyset \vdash \varphi$ belongs to the outer one: it
is a statement \emph{about} an internal derivation, made where we
stand, not a statement \oga{} makes.

Take the diagonal witness of \cref{thm:godel-fict}.  From outside,
we are not in any doubt about the underlying computation.  We know
its search hits at every argument --- at step zero, as
\cref{sec:metatheory:pt} observed --- so we know the totality
universal is $\omega$-true.  We also know \oga{} proves every
numeric instance of it, and does not prove its closure.  The
question ``does this search terminate everywhere?'' is settled,
and we are the ones who have settled it.

What the Fict covers is therefore not that question.  It covers
the absence of an internal certificate: the schematic derivation
that grounding the universal would demand, and that no
\oga-derivation supplies.  \oga{} does not learn the answer by
declaring the sentence decided; \ati{} confers status without
conferring either verdict.  \Cref{lem:totality-lem} makes the
point exactly --- the excluded middle over the sentence evaluates
to \ctrue{} while neither disjunct does.  The fiction asserts
\emph{that there is an answer}, never \emph{which}.

The phenomenon takes its sharpest form in \cref{thm:tier2}, which
exhibits a diagonal witness whose universal is $\omega$-true and
whose \emph{negation} may be adopted without contradiction.  A
fiction that remains compatible with the denial of a truth we
independently possess is not a fiction about the world.  It is a
fiction about certification.

One caution against generalizing.  It does not follow that every
gapped or undecided sentence is one whose answer we know from
outside.  We know it here only because the diagonal construction
hands us the halting fact alongside the unprovability.  The Liar
(\cref{sec:semantics:fict}) has no outer answer either, and
nothing above should be read as suggesting otherwise: the claim
concerns this family of witnesses, not ungroundedness in general.

\begin{rem}[no laundering; \texttt{omv\_iff\_transfer}]
\label{rem:no-laundering}
$\omega$-truth transfers across provable biconditionals: no
provable equivalence can smuggle a ceiling-bound sentence's content
into a provable variant.  Adopting a fiction therefore adds no new
grounded facts under any provable disguise --- adoption is
invisible at the grounded level.
\end{rem}

\subsection{Adopting the Anti-Fiction: the Safety Gradient Closes}
\label{sec:metatheory:tier2}

The fictional reading suggests a taxonomy of \emph{extensions}: an
added axiom may be $\omega$-true (safe --- tier 1), $\omega$-false
yet consistent (tier 2), or explosive (tier 3, the Liar's excluded
middle).  Tier 1 is inhabited by every fiction, per
\cref{rem:no-laundering} --- and its safety is not merely a remark
but a theorem suite, taken up in \cref{sec:metatheory:adopt}
below; tier 3 is inhabited by the classical paradoxes.  That
tier 2 is inhabited at all was open --- its classical witness,
adding $\neg G$ for a G\"odel sentence $G$, needs the deduction
theorem for the reductio, and grounded systems grant implication
introduction only for \emph{decided} antecedents.  In \rga{} the
diagonal sentence is not decided; the construction stalls.  In
\oga{} it is (\cref{thm:totality-decided}), and the construction
goes through --- decidedness unlocking implication introduction, the
mechanism flagged in \cref{sec:semantics:ladder}:

\begin{thm}[tier 2 is nonempty; \texttt{tier2\_nonempty}]
\label{thm:tier2}
There is a diagonal instance whose universal sentence $G$ is
$\omega$-true, while the extension of \oga{} by the hypothesis
$\neg G$ proves no contradiction.
\end{thm}

Adopting the \emph{fiction} $G$ is invisible at the grounded level;
adopting the \emph{anti-fiction} $\neg G$ is consistent but
$\omega$-unsound --- the two adoptions sit on opposite sides of the
gradient, and both are theorems about the same sentence.  The same
decidedness powers both constructions: decidedness giveth the
G\"odel-II-style argument and taketh away nothing here --- what it
may take away, namely self-verification, is the forward-looking
question of \cref{sec:concl}.

\subsection{Adopting the Fiction: Certified Stocks and Confluence}
\label{sec:metatheory:adopt}

The previous subsection closed the safety gradient's dangerous
side: adopting the \emph{anti}-fiction is consistent but
$\omega$-unsound.  This one arms the safe side, and in doing so
answers the operational questions that ``adoptable'' raises the
moment one takes it seriously as a practice rather than a status.
If fictions are to be adopted --- used as axioms of an extended
system --- then: what exactly warrants an adoption?  How many
fictions may one adopt at once?  And can two parties who adopt
independently ever collide?  Each answer is a theorem.

A word on reading the statements: an extension's derivability
relation --- \oga{} with a stock $E$ of adopted sentences --- is
literally derivability-from-hypotheses, so ``the extension proves
no contradiction'' and ``no sentence and its negation are both
derivable from $E$'' are the same statement read two ways.
Nothing beyond the existing proof system is involved: adoption
adds premises, never rules.

\begin{thm}[Adoption is safe;
  \texttt{omtrue\_hyp\_consistent},
  \texttt{omtrue\_fset\_consistent}]
\label{thm:adopt-safe}
Extending \oga{} by any well-formed $\omega$-true sentence ---
or by any finite stock of them --- yields no contradiction.
\end{thm}

The hypothesis is the honesty of the scheme, and it deserves the
emphasis: the warrant is $\omega$-truth, a \emph{metatheoretic}
condition, and adoption does not certify itself from inside ---
adopt falsehoods and no theorem protects you.  Call an adoption
\emph{certified} when its sentence carries this warrant (the
formalization names the profile \texttt{adoptable}: well-formed
and $\omega$-true), and \emph{uncertified} when the sentence is
adopted on its fictional status alone.  The distinction is not
pedantry, because status is symmetric: if $G$ is a fiction, so is
$\neg G$, and nothing about status distinguishes the true one.
The gradient's two sides are exactly this pair: adopting $G$ is
certified (tier 1, \cref{thm:adopt-safe}); adopting $\neg G$ is
uncertified --- singly consistent, as \cref{thm:tier2} showed,
but $\omega$-unsound.

\begin{thm}[Certified adoption is confluent;
  \texttt{adoption\_confluent}]
\label{thm:adopt-confluent}
If every member of the stock $E_1$ is certified adoptable, and
likewise every member of $E_2$, then the merged stock
$E_1 \cup E_2$ yields no contradiction.
\end{thm}

The proof is one appeal to \cref{thm:adopt-safe}, and that
brevity is the point.  The safety hypothesis is per-sentence,
never about a stock as a whole.  $\omega$-truth is truth in one
fixed semantics, and truth in a single model composes.  Certified
adoptions are therefore order-independent and parallelizable.  No
adoption ever needs its warrant upgraded to account for earlier
ones: every certified stock draws on the same single coherent
theory --- the $\omega$-true sentences --- and at most one of any
contradictory pair belongs to it.  Uncertified adoption enjoys no
such guarantee, and the gradient pair is the counterexample in
miniature.  $G$ and $\neg G$ are each singly consistent to adopt,
the first by \cref{thm:adopt-safe} and the second by
\cref{thm:tier2} --- both theorems of this paper --- yet their
union is inconsistent on its face.  \emph{Certified adoption is
confluent; uncertified adoption is not}, and it is the truth
anchor, not bookkeeping order, that makes the difference.

The subsection's results assemble into a single statement of the
paper's thesis:

\begin{thm}[The adoptable-fiction bundle;
  \texttt{gsent\_adoptable\_bundle}]
\label{thm:adopt-bundle}
Under the diagonal hypotheses of \cref{thm:godel-fict} --- the
search succeeds everywhere, the universal $G$ is unprovable ---
the sentence $G$ is simultaneously: certified decided;
irrefutable; unprovable; a genuine fiction; and safe to adopt,
the extension of \oga{} by $G$ proving no contradiction.
\end{thm}

This is incompleteness re-priced, in one theorem: the system
manufactures its own ideal elements, certifies their status,
computes their pedigree, and guarantees the safety of taking them
on.  What the working extension of \oga{} by adopted fictions
looks like at scale --- equalities adopted to recover
extensional reasoning, whole set-theoretic universes carried on
certified stocks --- is the companion set-theory paper's subject;
and where the warrant itself might live, one reflective level up,
is the closing question of \cref{sec:concl}.

\section{The Formalization}
\label{sec:formal}

Every result in this paper is machine-checked in
Isabelle/HOL~\cite{nipkow02isabelle}, extending the
grounded-deduction project's formalization session.  The
\oga-specific development comprises twenty theories, roughly
6{,}800 lines beyond the \rga{} base, with about 370 named
theorems and lemmas.  It contains no unproven obligations and adds
no axioms to HOL.  Each theorem header in this paper names the
corresponding Isabelle fact.

\subsection{Shape of the Development}
\label{sec:formal:shape}

Three structural decisions did the most work.
\emph{Case-aligned rule sets}: \oga's derivation relation restates
\rga's forty-two rules textually, plus \ati, so metatheoretic
inductions --- soundness, the certificate checker, the
strip/transport kit --- replay case for case, and the one new case
isolates exactly what \ati{} costs each invariant.
\emph{The one-step checker architecture}: proofs are judgment
lists checked entry-by-entry against their suffixes; the same
pattern instantiates for the proof checker and the semantic
certificate machinery, and (in the companion set-theory
development) twice more.
\emph{Delta presentation as a locale ladder}: the generic
quantifier-tier interfaces of \cref{sec:decidedness:tier} are
locales shared by all systems of the project, so the separation is
recorded as which tiers each system registers, not as ad-hoc
theorem pairs.

\subsection{Constructive Content}
\label{sec:formal:constructive}

Isabelle/HOL is classical, so the development inherits classical
logic whether or not it needs it, and a count of appeals to
contradiction would measure nothing.  The question worth asking is
where a theorem's \emph{content} genuinely lacks a computational
witness.  We audited the development on three markers: proof by
contradiction on a positive goal, choice or Zorn, and Hilbert
choice beyond decidable existence.

Non-constructivity turns out to be localized, and weaker than it
first appears.  Zorn's lemma occurs exactly twice, both inside the
completeness proof for the fafi semantics (\cref{thm:c1}), at the
maximal-consistent-extension step.  That appeal is not essential:
the term language is G\"odel-coded and therefore countable, and
countable Lindenbaum needs only a \emph{fixed} enumeration rather
than a chosen one, so the construction is choice-free in
principle.  Removing it is plausible and has not been done.  What
genuinely remains at that step is non-computability rather than
non-constructivity --- no effective procedure produces the
extension --- and either way it is unrelated to the diagonal
phenomena occupying the rest of the paper: it concerns extension,
not the distance between an r.e.\ set and a $\Pi_1$ one.
Soundness, the certificate checker, and the decidedness ladder use
none of the three markers.

What remained that looked classical turned out to be style, and it
has since been removed.  The diagonal results
(\cref{thm:incomplete,thm:godel-fict}) were originally proved by
contradiction from the non-enumerability of the halting complement
--- but the lemma they invoke is itself a \emph{constructive}
diagonal, building the offending point from the assumed
enumerator, and the contradiction wrapper then discarded it.
Recovering the witness turns the whole family into instances of a
single lemma:

\begin{quote}
Let $p$ be an r.e.\ property all of whose instances lie in the
complement of the halting diagonal.  Then $p$ misses a point of
that complement, namely the diagonal point of $p$'s own
enumerator.
\end{quote}

Each gap theorem supplies one soundness inclusion --- provability
implies membership in the complement --- and receives a named
witness in return.  All of them are now available in this form, on
both the \rga{} and the \oga{} side, and the witnessed versions
use a weaker hypothesis than the originals: the inclusion alone,
rather than coincidence of the two sets.  What looked like six
appeals to contradiction was one construction, repeated.

One residue is genuine rather than stylistic.  Even the witnessed
form is explicit only modulo a Hilbert choice in the indexing of
primitive recursive functions, so the gap point is concrete in
principle and opaque in practice.  That single shortcut is also
what keeps the compiler's code map from being primitive recursive,
and it is worth recording that the two costs are one cost: the
same convenience that makes an index unavailable as a computable
function of a function is what stops these witnesses from being
computed.

\subsection{What the Machine Contributed}
\label{sec:formal:lessons}

Three incidents from the development record show the proof
assistant acting as a research instrument rather than a notary,
and we report them because each changed the paper's content.

\emph{The admissibility repair.}  The fafi characterisation
(\cref{thm:c1}) was first attempted with a natural pointwise
notion of admissibility: a valuation should respect each
individually provable or refutable atom.  The machine refuted the
theorem in that form.  The failed proof exposed a valuation that
splits a pair of undecided atoms which the system, reasoning from
one as a hypothesis, provably entangles.  The corrected notion ---
coherent admissibility, a joint condition on the valuation's
commitments (\cref{sec:semantics:appx}) --- is not a repair of
the proof but a repair of the \emph{definition}, and it is the
one the paper's semantics now rests on.

\emph{The missing congruence.}  The head-step relation of
\cref{fig:rules-reflect} deliberately has no congruence rule
under binders: steps fire only at the head.  Several transport
arguments in the decidedness ladder appeared to need exactly such
a rule.  The formalization forced the honest alternative into the
open --- an explicit family of bridge lemmas that replay each
transport at the judgment level instead of rewriting under the
binder.  What looked like proof-engineering friction is in fact a
design invariant of the calculus, and \cref{sec:semantics:ladder}
states the ladder in the form the bridges made precise.

\emph{A semantics refuted before it was built.}  An earlier
design for the value semantics tracked, for each sentence, a
\emph{set} of possible classical truth values --- a
truth-functional scheme in the Kleene style.  Working the
soundness obligation for \ati{} on paper against that design
refuted it before any formalization began: possibility sets lose
the correlation between occurrences of the same unresolved
question, and the excluded-middle computation of
\cref{sec:semantics:fafi} is exactly what they cannot deliver.
The supervaluational design was adopted \emph{because} the
machine-bound obligation could not be met any other way --- the
classical distinction between truth-functional and
supervaluational gap theories, rediscovered operationally as a
proof obligation.

\subsection{AI Contribution Statement}
\label{sec:formal:ai}

Substantial parts of this work were carried out by an AI assistant
(Claude Fable from Anthropic), working under the direction and
review of the human author.  The AI contributions include the
mechanization of large parts of the \oga\ development --- the
system's rule layer and its soundness, the fact/fiction value
semantics and its coherence theorems, the decidedness and
tier-structure results, the recursive-enumerability metatheory with
its internalized proof checker, and the fictional-status
classifications --- as well as first drafts of this paper's text.
The human author directed the research, established the system
design and its revisions, posed the questions and thought
experiments that several of the developments answer, reviewed and
edited both proofs and prose throughout, and bears sole
responsibility for the paper's claims, framing, and treatment of
related work.  The collaboration is documented in the project
repository: working notebooks record the campaigns and design
decisions, and individual commits carry explicit AI co-author
attribution.  The formal results themselves are machine-checked by
Isabelle, whose kernel does not care who wrote the proof script.

\section{Related Work}
\label{sec:related}

\paragraph{Grounded and paracomplete truth.}
The grounding discipline descends from
Kripke's fixed-point theory of truth~\cite{kripke75outline};
\rga~\cite{rga-paper} made the quantifier clause reflective and
recursively enumerable.  \oga's fafi semantics applies the same
reflective move to \emph{supervaluation}: classical
supervaluational fixed points~\cite{kripke75outline,burgess86truth,
mcgee85truthlike} and related axiomatic systems such as
VF~\cite{cantini90theory} sit at $\Pi^1_1$ scale, whereas the
certificate-gated form stays r.e.  Axiomatic-truth costs and
comparisons follow
Halbach~\cite{halbach06axiomatizing,halbach18costs}.

\paragraph{Revenge.}
Paracomplete truth theories are characteristically vulnerable to
\emph{revenge}: the vocabulary needed to describe the gap
re-enacts the paradox, and the repair is an
ascent~\cite{beall08revenge} --- in Field's development, a
transfinite hierarchy of determinacy
operators~\cite{field08saving}.  \Cref{sec:metatheory:revenge}
argues that \oga{} occupies an intermediate position rather than
escaping the difficulty.  Its decidedness predicate is internal,
total on the reflective ground, and flat under iteration, so no
hierarchy is needed to ask whether a sentence is settled; but the
predicate that would carry Field's burden --- semantic valuedness,
the property whose failure \emph{is} the gap --- is deliberately
not internalized.  The retreat to a metalanguage survives; it has
been relocated to the one place the grounded discipline never
claimed to reach.

\paragraph{The $\omega$-rule and transfinite progressions.}
\ati{} reflects one $\omega$-fact --- decidedness --- into a
finitary rule, in contrast to the restricted
$\omega$-rule~\cite{shoenfield59omega} and to progressions along
transfinite iterations of reflection
principles~\cite{turing39ordinals,feferman62transfinite}.  The
comparison is instructive in both directions: progressions add
\emph{truth} at each stage; \ati{} adds only \emph{status}, which
is why the calculus stays r.e.\ and consistency is preserved by a
case-aligned replay.  Feferman's reflection principles assert, at
each stage, that the previous stage's theorems are \emph{sound} ---
each step is a leap of faith the new theory cannot itself justify.
Status-reflection makes a weaker leap: \ati{} asserts only that
surveyed questions \emph{have} answers, never which, and the
adoption discipline of \cref{sec:metatheory:adopt} then prices any
further leap sentence by sentence.

\paragraph{Self-verification and its limits.}
Willard's self-verifying systems~\cite{willard01self} obtain
consistency statements by weakening arithmetic; the grounded
family keeps arithmetic and weakens classical logic instead.  The
L\"ob obstruction~\cite{lob55solution} and its dependence on the
provability conditions frame the forward-looking question of
\cref{sec:concl}: at \oga{} the obstruction's known disabling
mechanism (undecided provability) is removed by \ati, with
consequences under active investigation.

\paragraph{Fictionalism and ideal elements.}
The pedigree apparatus gives formal content to a fictionalist
reading of ideal mathematical objects in the tradition of
Hilbert's ideal elements and modern
fictionalism~\cite{field08saving}: \oga's fictions are not
eliminable conveniences but certified, priced commitments ---
valued, jointly consistent with all grounded facts, and tracked
compositionally by the semantics itself.

\section{Conclusion and Outlook}
\label{sec:concl}

One rule separates \oga{} from \rga, and the separation is a tier
boundary: the reflective ground becomes decided, provability stays
recursively enumerable while $\omega$-truth deliberately does not,
and incompleteness returns transformed --- the G\"odel sentence is,
unconditionally, a genuine fiction with a computable pedigree.
The system does not merely fail to decide its diagonal sentences;
it certifies them as adoptable ideal elements and prices the
adoption.  All of this is machine-checked, and all of it flows
from a single premise-shape: decisions about infinite totalities
require surveying every instance, and \oga{} adopts exactly the
reflection of that survey into the object logic --- status,
never truth.

Three lines of work continue.  First, the \emph{graded set
theory} that \oga's certified searches support: sets with
provability-graded membership sustain the constructors of a
Zermelo--Fraenkel-style development with graded extensionality ---
where extensionality itself becomes an adoptable fiction with
computable pedigree --- and a Tarski-style universe constructor;
that development is the subject of the companion paper \emph{Fact
and Fiction in Grounded Set Theory} (in preparation), which
consumes this paper's semantics and adoption machinery as given.
Second, \emph{self-verification and the tier inversion}.  A
concurrent campaign investigates whether \rga{} proves its own
consistency; at \oga, mechanism-level analysis shows the classical
L\"ob obstruction's known disabling conditions are all lifted by
decidedness.  If both verdicts land as the evidence suggests, the
family exhibits a striking inversion: the weaker system
self-verifies while the stronger loses self-verification exactly
by gaining decidedness --- decidedness giveth the G\"odel-II-style
argument and taketh away consistency proofs.  The same decidedness
has already delivered the witness for the formerly-open
``unsound-yet-consistent'' tier of the safety gradient
(\cref{thm:tier2}): \oga{} plus the \emph{negation} of its
G\"odel sentence.  Working out the remaining consequences --- and
the derivability conditions a grounded provability predicate
actually satisfies --- is where this story goes next.

Third, the \emph{reflection rung}.  The adoption discipline of
\cref{sec:metatheory:tier2} takes its safety warrant from the
metatheory --- but the metatheory that supplies a warrant can
itself be a grounded system, reasoning about this one by
reflection, and there the warrant becomes a certified, generally
unprovable sentence: a fiction one level up, adoptable in turn.
Warrants stack into a tower, and the rule that would consume them
--- from a derivation of the warrant, license the extended system
--- is structurally the same move that builds \oga{} from \rga:
\ati{} consumes instance-wise decidedness certificates, the
warrant rule would consume instance-wise provability certificates.
This suggests the GA family is not four systems but the first
rungs of a generated progression, in the lineage of Turing's
ordinal logics and Feferman's transfinite
progressions~\cite{turing39ordinals,feferman62transfinite}, with
the paracomplete difference that every rung's leap of faith is a
certified, pedigreed object rather than a bare axiom.  Making that
rule precise is the family's natural next system.

\paragraph{Certification, and what grounding does not buy.}
\label{sec:concl:pcc}
The apparatus assembled here has an obvious affinity with
proof-carrying code~\cite{necula97proof}.  \oga{} proofs are
numbers; the checker is primitive recursive (\cref{thm:re}); and
$\omega$-soundness says an accepted certificate entails its
sentence.  Gate a computation on the check and the gated
computation is unconditionally safe.  For termination
specifically, \cref{thm:pt} gives the corresponding guarantee
directly: a step function whose totality \oga{} affirms really
does halt everywhere.

We think it is worth stating plainly how little of this is ours.
The pattern --- decidable checker, plus soundness, therefore safe
gate --- is what proof-carrying code has always been, and it
requires nothing of a logic beyond being sound and having a
decidable proof relation.  Every system in the tradition supplies
it.  Neither grounding nor \ati{} is doing work here, and a
version of this paragraph could be written about any sound
formalism.

Two further deflations are worth being explicit about, since the
reflective machinery invites overclaiming.  First, the gate is the
\emph{external} checker, not the internal provability predicate:
\cref{thm:prov-char} is a meta-theorem, and internally \oga{} is
no better placed than \pa{} --- it does not prove its own
soundness schema, \cref{thm:prov-not-internal} blocks one half
outright, and the constraint noted in \cref{sec:metatheory:re}
blocks the rest short of a consistency proof.  Reflection does not
make the gate self-certifying.  Second, the interesting comparison
is not with classical logic but with constructive type
theory~\cite{huet16coq,bove09brief,chlipala13certified}, where
extraction already yields certificates carrying computational
content.  Against that baseline we can identify no advantage in
what is proved here.

What might genuinely differ, offered as conjecture rather than
result, is a matter of trusted base and of distance.  \oga's
certificates are numerals and its checker is a primitive recursive
function, which is a smaller and more uniform object to trust than
the kernel of a dependent type theory --- though that is an
argument for small kernels, made before and independently of
grounding.  And because $\omega$-truth is defined by the
evaluation relation rather than axiomatized alongside it, ``$\varphi$
is true'' already means ``the process grounds $\varphi$'', which
shortens the distance between the logic and the operational model
a certificate is ultimately about.  But \oga's operational model
is its own term language, not a machine; compiling real machine
semantics into it would reopen exactly the gap that foundational
approaches to proof-carrying code were designed to close.

The honest summary is that the certification story belongs to the
grounded program rather than to this paper's rule.  Every
ingredient of it --- soundness, the r.e.\ checker, sound totality
affirmation --- has an \rga{} analogue, and \ati{} is not what
delivers any of them.  \ati{} bears on certification in exactly
two places, both modest.  Whether it enlarges the class of
programs certifiable as terminating is the strictness question of
\cref{sec:metatheory:pt}, which is open.  And
\cref{cor:decd-decided} lets the system reason internally by cases
on whether a certificate exists, even where it cannot say which
case obtains --- a reasoning device whose practical value we have
not tested.  A serious assessment of grounded arithmetic as a
certification substrate would be a different paper, and would have
to be argued against constructive type theory on its merits rather
than against classical logic on ours.

\bibliographystyle{plain}
\bibliography{logic}

\appendix
\section{The Fafi Machinery, and the Proof of the Characterisation}
\label{sec:semantics:appx}

This appendix sketches the machinery behind
\cref{thm:c1} --- the coincidence of provability with super-truth
--- at the level of ideas.  The full development is in the
formalization (theories \texttt{OGA\_Fafi\_Atoms},
\texttt{OGA\_Fafi\_BA}, \texttt{OGA\_Fafi\_Eval}, and the
\texttt{OGA\_Fafi\_C1} chain), and every claim below names its
machine-checked counterpart.

\paragraph{Evaluation.}  The evaluation relation \texttt{feval}
assigns values by induction on syntactic shape, with one clause per
form.  Grounded atomic sentences evaluate to constants.  The
connective clauses combine values through the Boolean algebra.  The
one distinctive clause is the gate for universals described in
\cref{sec:semantics:fafi}: a universal whose body is certified
decided may take its own named atom as a value.  Two structural
points matter later.  Evaluation is a \emph{relation}, not a
function --- a sentence may carry several values, and coherence
between them is a theorem rather than a stipulation.  And the atom
discipline is stable under substitution
(\texttt{OGA\_Fafi\_Atoms}), which is what lets valued sentences
pass under binders.

\paragraph{Coherent admissibility.}  A valuation $w$ assigns a
truth value to every atom.  Call a finite hypothesis context
\emph{honest} for $w$ (\texttt{honest}) when each of its members is
a signed decided atom that $w$ asserts: either a certified
universal that $w$ makes true, or the negation of one that $w$
makes false.  The valuation is \emph{admissible} (\texttt{admiss})
when no honest context is contradictory --- from hypotheses $w$
itself endorses, the system can never derive both a sentence and
its negation.  This is the exact content of the gloss in
\cref{thm:c1}: an admissible valuation never contradicts what the
system can actually derive, even hypothetically.

The definition earns a design note, because it is the second
attempt, and the machine refuted the first.  The naive notion ---
admissibility as respecting each \emph{individually} provable or
refutable atom --- makes the characterisation \emph{false}: a
valuation could split a pair of undecided atoms that the system,
reasoning from one as hypothesis, provably entangles.  The
formalization found the counterexample during the campaign, and
the corrected, context-based notion above is what survives it.
Coherence is a joint condition on the valuation's commitments, not
a pointwise one.

\paragraph{Soundness: the truth lemma.}  The soundness half
(\texttt{c1\_soundness}) shows a provable sentence is super-true.
It runs through a \emph{truth lemma} (\texttt{TL},
strengthened to an equivalence in \texttt{TL\_iff}): for an
admissible valuation, the value a sentence takes agrees with
provability from honest contexts.  The proof is structural, and
notably frugal.  It uses no disjunction property and no maximal
theories: a total valuation already decides every atom, the
connective clauses lift that completeness through the syntax, and
the primeness that classical arguments buy with maximality comes
free from disjunctive syllogism.

\paragraph{Completeness: Zorn, without the deduction theorem.}
The completeness half (\texttt{c1\_completeness}) shows a
super-true sentence is provable, and its interest lies in what it
must avoid.  The classical route --- extend $\{\neg\varphi\}$ to a
maximal consistent theory --- is unavailable twice over.  Grounded
logics grant implication introduction only for decided antecedents
(\cref{sec:bg}), so the usual deduction-theorem bookkeeping fails.
And a sentence's negation may itself be ungrounded, so the seed
set $\{\neg\varphi\}$ is not a legitimate starting point.  The
formalization instead runs Zorn's lemma directly on the family of
atomic sets that are \emph{consistent and do not prove} $\varphi$
--- no negation seed at all.  The extension step adds a decided
atom by case analysis on its excluded middle, which is exactly
what decidedness licenses, in place of the unavailable deduction
theorem.  A maximal such set is complete on decided atoms
(\texttt{both\_prove\_collapse}), so it induces a valuation.  That
valuation is admissible by construction, and it refutes $\varphi$
--- contradicting super-truth.  Coherent models of consistent
atomic sets exist by the same construction (\texttt{lindenbaum},
\texttt{coherent\_model\_exists}).

\paragraph{The bridge to the $\omega$-completion.}  The two
semantics of \cref{sec:semantics} meet in one valuation.  Define
$w_\omega$ to assert exactly the atoms whose sentences are
$\omega$-true (\texttt{womega}).  This valuation is admissible,
and a value read at $w_\omega$ agrees with $\omega$-truth of its
formula.  Super-truth therefore implies $\omega$-truth: whatever
is true under \emph{every} admissible resolution of the open
questions is in particular true under the resolution the
$\omega$-model supplies.  The fafi semantics refines the
$\omega$-completion rather than competing with it.

\end{document}